\documentclass[12pt]{report}
 
\usepackage{amsfonts} 
\usepackage{amsbsy}

\newcommand{\eh}{\hspace{.05in}}  
\newcommand{\ih}{\'\i}  
\newcommand{\C}{\mathbb{C}}
\newcommand{\E}{\mathbb{E}}
\newcommand{\N}{\mathbb{N}}
\newcommand{\R}{\mathbb{R}}
\newcommand{\Z}{\mathbb{Z}}
\newcommand{\A}{\mathcal A}
\newcommand{\cC}{\mathcal C}

\newcommand{\F}{\mathcal F}
\newcommand{\K}{\mathcal K}
\newcommand{\cR}{\mathcal R}
\newcommand{\cT}{\mathcal T}
\newcommand{\cZ}{\mathcal Z}
\newcommand{\sD}{\mathsf D}
\newcommand{\sK}{\mathsf K}
\newcommand{\m}{\frac{1}{2}} 
\newcommand{\dd}{\partial}  
\newcommand{\ds}{\displaystyle}  
\newcommand{\BE}{\begin{equation}}  
\newcommand{\EE}{\end{equation}}  
\newcommand{\af}{\alpha}   
\newcommand{\bt}{\beta}
\newcommand{\eps}{\varepsilon}
\newcommand{\ka}{\kappa}
\newcommand{\ld}{\lambda}
\newcommand{\Ld}{\Lambda}
\newcommand{\sg}{\sigma}
\newcommand{\Sg}{\Sigma}
\newcommand{\bz}{\bar{z}}
\newcommand{\bg}{\bar{g}}
\newcommand{\bS}{\bar{S}}
\newcommand{\ta}{\mathtt a}

\newcommand{\bba}{\boldsymbol a}
\newcommand{\bbg}{\boldsymbol g}
\newcommand{\bbf}{\boldsymbol\phi}
\newcommand{\bbx}{\boldsymbol x}
\newcommand{\bbz}{\boldsymbol z}
\newcommand{\bbsg}{\boldsymbol\sg}
\newcommand{\bbK}{\boldsymbol K}
\newcommand{\Lim}[1]{\lower5.5pt\hbox{${{\ds\lim}\atop^{#1}}$} \ } 
\newcommand{\fns}{\footnotesize}

\begin{document}
\ \\   
{\large{\bf Triply periodic minimal surfaces which converge to the}}

\centerline{\large{\bf Hoffman-Wohlgemuth example}}
\ \\ 
\centerline{P{\fns LINIO} S{\fns IM\~OES} \& V{\fns AL\'ERIO} R{\fns AMOS} B{\fns ATISTA}}
\ \\
{\bf Abstract.} We get a continuous one-parameter new family of embedded minimal surfaces, of which the period problems are two-dimensional. Moreover, one proves that it has Scherk's second surface and Hoffman-Wohlgemuth's example as limit-members.
\\
\\
{\bf 1. Introduction} 
\\  

A continuous family $\F$ of complete embedded minimal surfaces can play an important role in the development of their global theory. One of the most beautiful examples is {\it the genus one helicoid}, of which embeddedness was proved in 2000 by Hoffman, Wolf and Weber (see [HMM]), {\it seven} years after its discovery by Hoffman, Karcher and Wei (see [HKW]). Weber first showed that it was a limit-member of such an $\F$, and then used [HKW], [HPR] and the {\it maximum principle} to conclude his proof. With that, he finally added a second example of complete minimal submanifold of $\R^3$ {\it with only one end}, besides the helicoid. To date, one has not found any further examples of this kind yet.

Sometimes, one can find $\F$ enclosing all surfaces of a certain class. For instance, in 2005 P\'erez, Rodr\ih guez and Traizet proved that any doubly periodic minimal torus with parallel ends is an interior point of a cube $\F$. Its edges stand for either Scherk's or Riemann's examples, while each vertex is either the catenoid or the helicoid (see [PRT]). Such families are essential to understand the moduli space of minimal surfaces. Roughly saying, in the same connected component of this space, any two surfaces can be continuously deformed, one into the other and always keeping the minimality condition. 

At this point, we remark that the above references deal with \it two-dimen- sional \rm period problems. By this concept we do not count L\'opez-Ros parameters, and that dimension has been the highest in which one succeeds in finding a non-trivial {\it explicit} $\F$. To date, there still remain only few such examples, while many $\F$'s were obtained from {\it one-dimensional} period problems (see [HK], [K1-3] and [V1]).  

By the way, [V1] builds a strong parallel to this present work, for there one proves that Scherk's second surface and Callahan-Hoffman-Meeks' [C] are limit-members of a {\it unique} $\F$, in the sense that it encloses all the examples presented therein. In this paper we show that handle addition is possible for that {\it whole} $\F$, with one limit-member being an example from Hoffman and Wohlgemuth (see [HW] and [SV]).

If one seeks after a new isolated surface with less than three period problems, then handle addition is an old and widely known technique, though not always successful. In this work, however, we not only present a full study of a continuous family of new surfaces, but also do it practically {\it without} computations. Instead, geometric arguments are intensively used, many of them profiting from former results like [MR] and [V1]. By studying periods, one takes homotopic curves based on a {\it best-choice} procedure, detailed in Section 6.  
\input epsf  
\begin{figure} [ht]  
\centerline{  
\epsfxsize 9cm  
\epsfbox{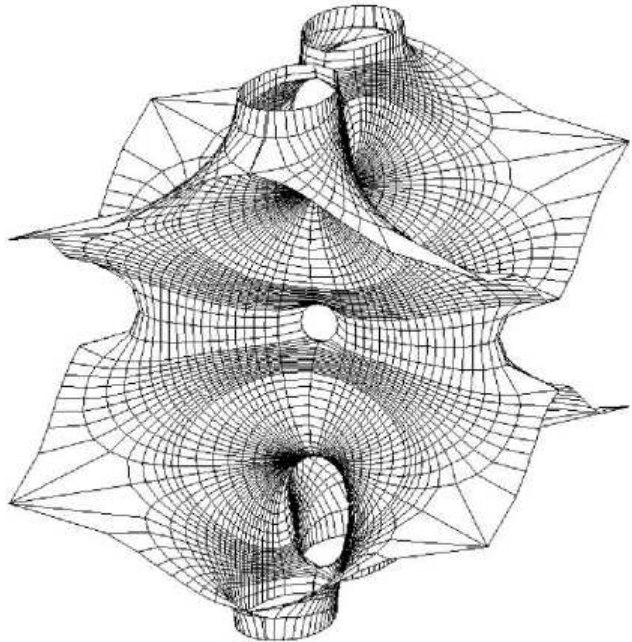}}
\caption{Fundamental piece of a triply periodic Costa surface with handles.}  
\end{figure} 
 
Let us first consider Figure 1. The main goal of this paper is then to prove the following:  
\\ 
\\ 
{\bf Theorem 1.1.} \it There exists a one-parameter family of complete triply periodic minimal surfaces in $\R^3$ such that, for any member of this family the following holds:
\\
(a) The quotient by its translation group $G$ has genus 7.
\\
(b) The whole surface is generated by a fundamental piece, which is a surface with boundary in $\R^3$. The boundary consists of eight planar curves of vertical reflectional symmetry and four planar curves of horizontal reflectional symmetry. The fundamental piece has a symmetry group generated by two vertical planes of reflectional symmetry and two line segments of 180$^\circ$-rotational symmetry.  
\\
(c) By successive reflections in the boundary of the fundamental piece one obtains the triply periodic surface.
\\
(d) All members in the family are embedded in $\R^3$. Moreover, it has two limit-members: the Hoffman-Wohlgemuth example of genus 5 and two side-by-side copies of Scherk's doubly periodic surface.\rm
\\

This work was supported by FAPESP grant numbers 00/07090-5, 01/05845-1 and 05/00026-3.
\\
\\
{\bf 2. Preliminaries}
\\

In this section we state some basic definitions and theorems. Throughout this work, surfaces are considered connected and regular. Details can be found in [K3], [LM], [N] and [O].
\\
\\
{\bf Theorem 2.1.} \it Let $X:R\to\E$ be a complete isometric immersion of a Riemannian surface $R$ into a three-dimensional complete flat space $\E$. If $X$ is minimal and the total Gaussian curvature $\int_R K dA$ is finite, then $R$ is biholomorphic to a compact Riemann surface $\overline{R}$ punched at a finite number of points.\rm
\\
\\
{\bf Theorem 2.2.} (Weierstrass representation). \it Let $R$ be a Riemann surface, $g$ and $dh$ meromorphic function and 1-differential form on $R$, such that the zeros of $dh$ coincide with the poles and zeros of $g$. Suppose that $X:R\to\E$, given by
\[
   X(p):=Re\int^p(\phi_1,\phi_2,\phi_3),\eh\eh where\eh\eh
   (\phi_1,\phi_2,\phi_3):=\m(g^{-1}-g,ig^{-1}+ig,2)dh,
\]
is well-defined. Then $X$ is a conformal minimal immersion. Conversely, every conformal minimal immersion $X:R\to\E$ can be expressed as above for some meromorphic function $g$ and 1-form $dh$.\rm
\\
\\
{\bf Definition 2.1.} The pair $(g,dh)$ is the \it Weierstrass data \rm and $\phi_1$, $\phi_2$, $\phi_3$ are the \it Weierstrass forms \rm on $R$ of the minimal immersion $X:R\to X(R)\subset\E$.\\
\\
{\bf Theorem 2.3.} \it Under the hypotheses of Theorems 2.1 and 2.2, the Weierstrass data $(g,dh)$ extend meromorphically on $\overline{R}$.\rm
\\

The function $g$ is the stereographic projection of the Gau\ss \ map $N:R\to S^2$ of the minimal immersion $X$. It is a covering map of $\hat\C$ and $\int_SKdA=-4\pi$deg$(g)$. These facts will be largely used throughout this work.
\\ 
\\  
{\bf 3. The symmetries of the surface and the elliptic $Z$-function}  
\\   
 
Let us consider Figure 1, which represents the fundamental piece of a triply periodic surface $S$. If $G$ denotes its translation group, then $S/G$ is a compact Riemann surface of genus 7 that we call $\bS$ (see Figure 2(a)). Let $\rho$ be the map from $\bS$ to its quotient by 180$^\circ$-rotation around the $x_3$-axis. Then, the Euler-Poincar\'e characteristic of $\rho(\bS)$ is given by $\chi(\rho(\bS))=\frac{\ds\chi(\bS)}{\ds 2}+6=0$. Because of this, $\rho(\bS)$ is a torus that we call $T$. This torus must be rectangular because of the following argument. The horizontal reflectional symmetries of $\bS$ are inherited by $T$ through $\rho$, and there are two curves which remain invariant under any of these symmetries. Then, the fixed-point set has two components and this only happens for the rectangular torus. 
 
\input epsf  
\begin{figure} [ht] 
\centerline{ 
\epsfxsize 12cm 
\epsfbox{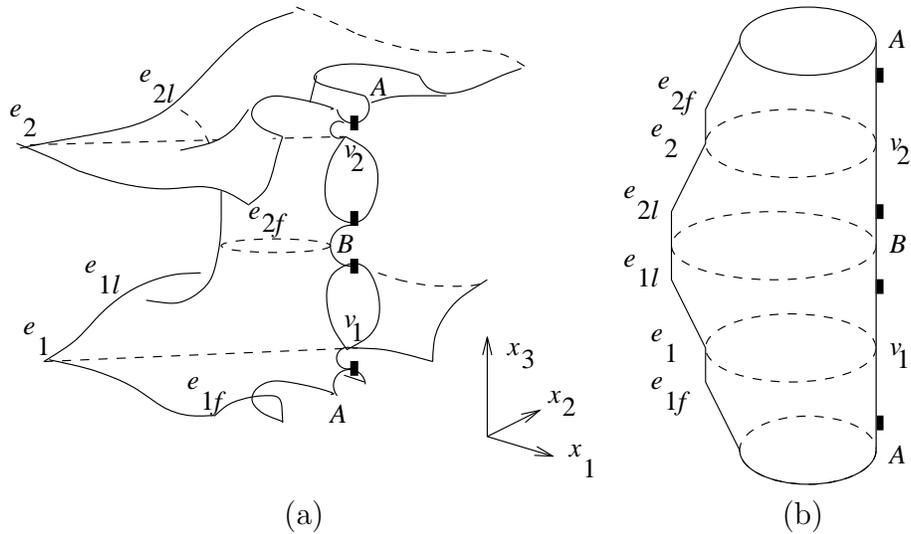}} 
\hspace{1.8in}(a)\hspace{2.4in}(b)
\caption{(a) Half of $\bS$; (b) the torus $T$.}  
\end{figure} 
 
The surface $\bS$ has two other 180$^\circ$-rotational symmetries, namely the ones around the $x_1$- and $x_2$-axes. The torus $T$ has these two symmetries as well. Let $r$ be the 180$^\circ$-rotational symmetry around the $x_1$-axis. The quotient of $T$ by $r$ is conformally $S^2$. After we fix an identification of $S^2$ with $\hat\C$, we finally obtain an elliptic function $Z:T\to S^2$. 
 
Consider Figure 2(b) and the points of the torus $T$ represented there. These correspond to special points of $\bS$, indicated in Figure 2(a) (they were given the same names). Let $Z:T\to S^2$ be the elliptic function with $Z(e_1)=1/Z(e_2)=0$ and $Z(v_1)=1/Z(v_2)=a$, where $a$ is a real value in $(0,1)$ (these functions coincide with $\cos\af\cdot\wp+\sin\af$ described in [K, p.40]).
   
\input epsf  
\begin{figure} [ht]  
\centerline{  
\epsfxsize 10cm  
\epsfbox{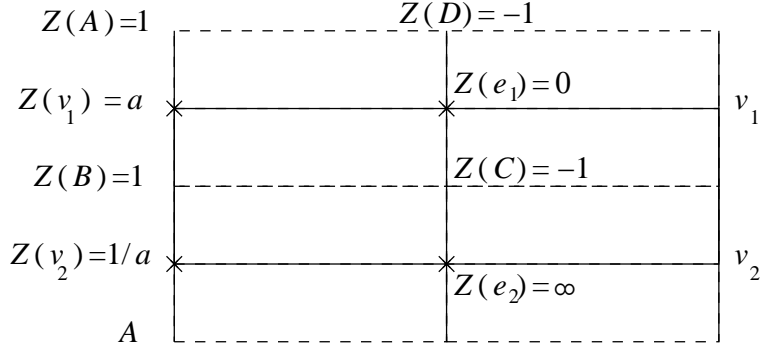}}  
\caption{The torus $T$ with values of $Z$ at special points on it.}  
\end{figure} 
  
Now we summarise some important properties of the function $Z$ (see Figure 3). It is real on the bold lines (and nowhere else), and $|Z|=1$ on the dashed lines (and nowhere else). It has exactly four branch points, marked with $\times$ in Figure 3. At the points $A$ and $B$ the function $Z$ takes the value $1$ and at the points $C$ (centre) and $D$, the value $-1$. 
\\
\\ 
{\bf 4. The $z$-function on $\bS$ and the Gauss map in terms of $z$} 
\\ 
 
In this section we start by studying the necessary conditions for the existence of a minimal surface like in Figure 1. They will lead to an algebraic equation for the compact Riemann surface $\bS$, together with Weierstrass data on it. From this point on, our problem will be concrete. We shall have to prove that the algebraic equation really corresponds to $\bS$ in terms of its genus and symmetries. Afterwards, we shall have to prove that the Weierstrass data really lead to a minimal embedding of $\bS$ in $\R^3/G$ with the expected properties: symmetry curves, periodicity, etc.  
 
Let us call $S$ the surface represented in Figure 1 and suppose that it is a minimal immersion of $\bS$ in $\R^3/G$. In this case, we make use of the previous section and consider the functions $\rho:\bS\to T$ and $Z:T\to\C$. Let us define $z:=Z\circ\rho$. Both functions $Z$ and $\rho$ have degree 2, then $z$ is a function on $\bS$ of degree 4 (see Figure 4(a)). In this picture one sees that $z$ takes on special values $b\in(a,1)$ and $-x\in(-1,0)$ on $\bS$.

\input epsf  
\begin{figure} [ht]  
\centerline{  
\epsfxsize 16cm  
\epsfbox{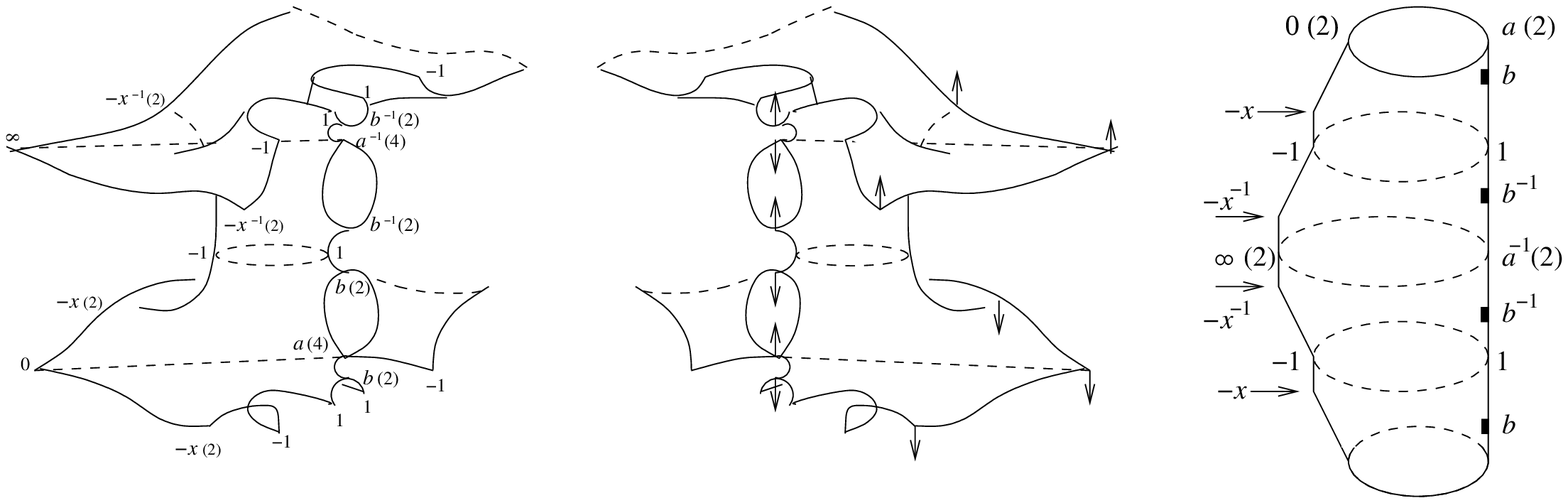}}  
\hspace{1in}(a)\hspace{1.5in}(b)\hspace{2in}(c)  
\caption{(a) Values of $z$ at special points; (b) The corresponding normal}    
\hspace{.8in}vector at these points; (c) the corresponding values of $Z$ on $T$. 
\end{figure}   
 
We are supposing that $S$ is a minimal immersion of $\bS$ in $\R^3/G$. In this case, the Gauss map on $S$ must lead to a meromorphic function $g$ on $\bS$, as Figure 4(b) suggests. We are going to define {\it multiplicity} as the branch order plus one. Then, the expected correspondence between the values of $z$ and $g$ (including their multiplicities) is indicated in Figure 4(a) and 4(b). Therefore, one settles the following relation:  
\BE    
   g^4=z\biggl(\frac{1-az}{z-a}\biggl)\biggl(\frac{b-z}{bz-1}\biggl)^2\biggl(\frac{z+x}{xz+1}\biggl)^2.  
\EE 

From now on we define $\bS$ as a general member of the family of compact Riemann surfaces given by (1). These surfaces have genus 7, because of the following argument: each value $z\in\{a^{\pm 1},0^{\pm 1}\}$ represents 1 branch point of multiplicity 4 on $\bS$, and each value $z\in\{-x^{\pm 1},b^{\pm 1}\}$ represents 2 different branch points of multiplicity 2 on $\bS$. This function is a four-sheet branched covering of the sphere. Therefore, by the Riemann-Hurwitz formula, the genus of $\bS$ is   
\[               
   \frac{4\cdot 1\cdot(4-1)+4\cdot 2\cdot(2-1)}{2}-4+1=7.    
\] 
 
Some involutions of $\bS$ are summarised in Table (2). This table includes the differential $dh$ which will be discussed in the next section. 
\BE 
\begin{tabular}{|c|c|c|c|c|}\hline 
 $   $&$ {\rm involution}       $&$ z{\rm -values} $&$ g\in$&$ dh(\dot{z})\in$\\\hline\hline
 $ 1 $&$ (z,g)\to(\bz,\bg)      $&$ -1<z<-x        $&$ \R  $&$ \R            $\\ \hline  
 $ 2 $&$ (z,g)\to(\bz,-\bg)     $&$ -x<z<0         $&$ i\R $&$ \R            $\\ \hline   
 $ 3 $&$ (z,g)\to(\bz,\pm i\bg) $&$ 0<z<a          $&$ \pm\sqrt{i}\R $&$ i\R $\\ \hline   
 $ 4 $&$ (z,g)\to(\bz,-\bg)     $&$ a<z<b          $&$ i\R $&$ \R            $\\ \hline    
 $ 5 $&$ (z,g)\to(\bz,\bg)      $&$ b<z<1          $&$ \R  $&$ \R            $\\ \hline    
 $ 6 $&$ (z,g)\to(1/\bz,1/\bg)  $&$ z\in S^1       $&$ S^1 $&$ i\R           $\\ \hline 
\end{tabular} 
\EE   
 
We have just proved that the values of $g$ on all special curves of $\bS$ are consistent with the expected unitary normal on the minimal surface $S$ in $\R^3/G$. 
\\ 
\\ 
{\bf 5. The height differential $dh$ in terms of $z$}  
\\ 
 
Now we need an expression for the differential form $dh$. The surface has no ends and because of this $dh$ is holomorphic. Its zeros are exactly the ones where $g=0$ or $g=\infty$ and all have multiplicity 1 (i.e., branch order 0). If we consider the differential form $dz$, then it would be sufficient to divide it by a function on the surface with double zeros at $z\in\{0,a^{\pm 1}\}$ and a pole of multiplicity 6 at $z=\infty$. This function will turn out to be the pull-back by $\rho$ of another function, that we call $V$, on the torus $T$.     
 
Since $0^{\pm 1}$ and $a^{\pm 1}$ are the only branch values of $Z$, all of them of order one, then the torus $T$ can be algebraically described by the equation 
\BE
   V^2=Z(Z-a)(Z-1/a). 
\EE
Now, $V\circ\rho$ has exactly the zeros and poles on $\bS$ with the expected multiplicities. We can take $v:=V\circ\rho$. This means that $v$ is a well-defined square root of the function $z(z-a)(z-1/a)$ on $\bS$. For instance, $v/z=:\sqrt{z+1/z-a-1/a}$.

Finally, we need to establish a proportional constant to determine $dh$ by means of $dz/v$. On the straight lines of the surface, where $0<z^{\pm 1}<a$, the coordinate \eh $x_3=Re\int dh$ \eh must be constant. Then $Re\{dh\}$ is zero there. Because of this we choose the proportional constant to be $i$, namely  
\BE 
   dh=\frac{idz}{v}=\frac{idz/z}{\sqrt{z+1/z-a-1/a}}.
\EE 
 
At this point we have reached concrete Weierstrass data $(g,dh)$ on $\bS$, defined by (1) and (4), with $x$, $a$ and $b$ satisfying the following inequalities
\BE
   0<a<b<1\eh\eh{\rm and}\eh\eh 0<x<1.
\EE

Now our task will be the demonstration of the following: let $S$ be the minimal immersion of $\bS$ given by these Weierstrass data. Then $S$ leads to the expected surface of which the fundamental piece is represented in Figure 1. In other words, we need to show that $S$ really has all the symmetry curves and lines of our initial assumptions, and the fundamental piece of $S$ has {\it no periods}, as indicated in Figure 1. This second task will be discussed in the next section. Now we analyse the symmetries of $S$. 
 
From (1) and (4) we see that all the $z$-curves listed in (2) are geodesics, because $g(z)$ is contained either in a meridian or in the equator of $S^2$, {\it and} $dh(\dot{z})$ is contained in a meridian of $S^2$. Moreover, the geodesics are straight lines if $0<z^{\pm 1}<a$, because in this case $\frac{\ds dg(\dot z)}{\ds g(z)}\cdot dh(\dot z)\in i\R$. Otherwise we shall have $\frac{\ds dg(\dot z)}{\ds g(z)}\cdot dh(\dot z)\in\R$ and the corresponding geodesics will be planar. Therefore, $S$ has all the expected symmetries. 
\\ 
\\
{\bf 6. Solution of the period problems} 
\\ 
 
The triply periodic minimal surface $S$ is generated by its translation group $G$ applied to a fundamental piece. Its right half is shown in Figure 5(a). The fundamental domain for the full symmetry group of the minimal surface is the shaded region represented on Figure 5(a). 
 
\input epsf 
\begin{figure} [ht] 
\centerline{ 
\epsfxsize 13.5cm 
\epsfbox{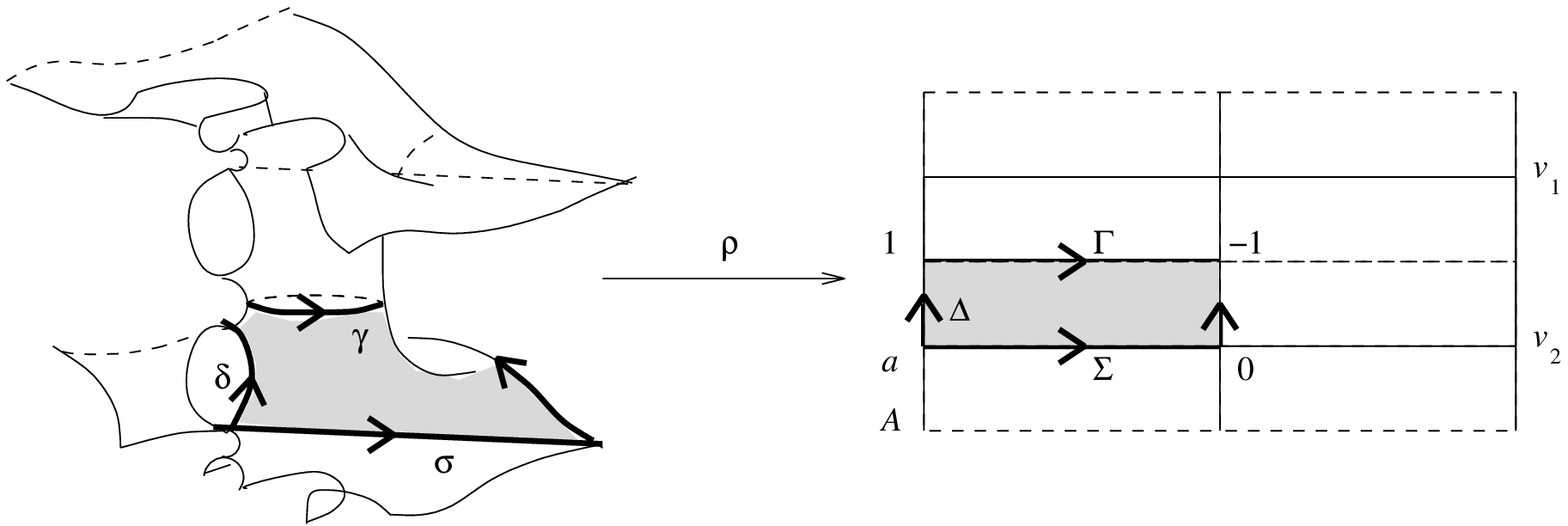}} 
\hspace{1in}(a)\hspace{3in}(b) 
\caption{(a) The right half of the fundamental piece;} 
\hspace{1.55in}(b) Its corresponding image under $\rho$. 
\end{figure} 
 
Since $S$ has no ends, we just need to analyse the period vector given by $Re\oint(\phi_1,\phi_2,\phi_3)$ on the curves of the homology of $\bS$. This task is very similar to the analysis done in [V1, p.80-81] and will be skipped here. We conclude that just two period problems remain to be solved, namely
\BE
   Re\int_{\gamma}\phi_2=0\eh\eh{\rm and}\eh\eh Re\int_{\delta}\phi_2=0,
\EE
where $\gamma$ and $\delta$ are represented in Figure 5(a). The branches of the square root need to be chosen in accordance with Figures 5(a) and 5(b). This choice is indicated in Figure 6. 

\input epsf 
\begin{figure} [ht] 
\centerline{ 
\epsfxsize 14cm 
\epsfbox{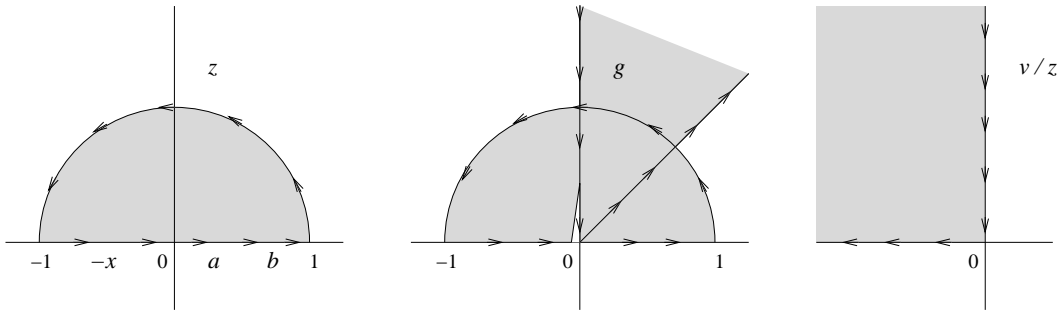}} 
\caption{The images of $|z|<1<1+Im\{z\}$ under $g$ and $v/z$.} 
\end{figure} 

The curve $\gamma$ can be explicitly given by $z(t)=z\circ\gamma(t)=e^{it},0<t<\pi$. If we define $\Gamma:=\rho\circ\gamma$, then $Z\circ\Gamma(t)=z\circ\gamma(t)$. We establish the 4$^{\rm th}$-root on $z(t)$ of each factor in (1) as indicated in Figure 7.
 
\input epsf 
\begin{figure} [ht] 
\centerline{ 
\epsfxsize 13cm 
\epsfbox{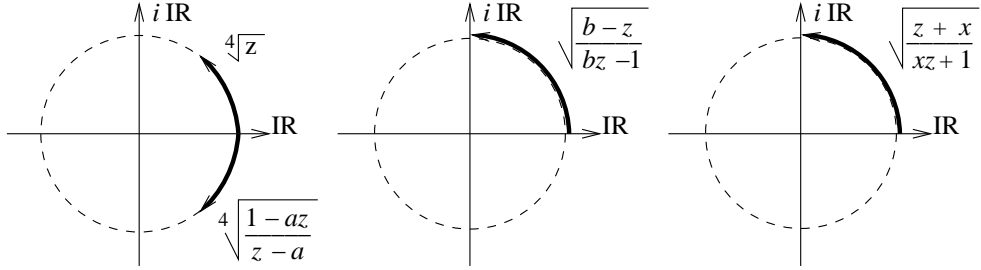}} 
\caption{The 4th-roots on $z(t)$ of the factors in (1).} 
\end{figure} 
 
The condition $Re\int_{\gamma}\phi_2=0$ will then be equivalent to 
\BE
   \m\int_{0}^{\pi}(g+g^{-1})|dh|=\int_{0}^{\pi}Re(g)|dh|=0,
\EE
where 
\BE
   dh=\frac{idt}{\sqrt{a+1/a-2\cos t}}\eh\eh{\rm and}\eh\eh
   g=g(z(t))\in S^1.
\EE

It is not difficult to see that $Re(g(t))$ is increasing with $x$ and decreasing with $b$. Let us now vary $b$ in the interval $(a,1)$. From {\it Lebegue's dominated convergence theorem}, at the extremes we have the following equalities for $I_\gamma:=\int_{0}^{\pi}Re(g)|dh|$:
\BE
   I_\gamma|_{b=a}=\int_0^{\pi}Re\biggl\{\sqrt[4]{z}\cdot
   \sqrt[4]{\frac{z-a}{1-az}}\cdot\sqrt{\frac{z+x}{xz+1}}\biggl\}
   \frac{dt}{\sqrt{a+1/a-2\cos t}} 
\EE
and
\BE
   I_\gamma|_{b=1}=-\int_{0}^{\pi}Im\biggl\{\sqrt[4]{z}\cdot
   \sqrt[4]{\frac{1-az}{z-a}}\cdot\sqrt{\frac{z+x}{xz+1}}\biggl\}
   \frac{dt}{\sqrt{a+1/a-2\cos t}}.
\EE

Both functions in (9) and (10) are still increasing with $x$. Let us analyse the integrand of (10). It is easy to prove that
\BE
   \tan Arg\biggl\{z\cdot\frac{1-az}{z-a}\biggl\}=
   \frac{2\sin t\cdot(a\cos t-1)}{a+1/a}.
\EE 

Hence, at $x=1$ the integrand of (10) will be always positive and consequently $I_\gamma|_{(b,x)=(1,1)}>0$ for any $a\in(0,1)$. It is not difficult to see that $a^{-1/2}I_\gamma|_{(b,x)=(1,0)}$ is negative for $a$ close to zero, while it diverges to $+\infty$ when $a$ approaches $1$. Notice that the factor $Im\{\cdot\}$ is monotonely decreasing with $a$. For any fixed $a\in(0,1)$, it changes sign at a certain unique $t_a\in(0,\pi)$. Now consider a value $a=\af$ where $I_\gamma|_{(b,x)=(1,0)}$ vanishes. If one takes $p:=(1/\af-\af)/(1/\af+\af-2\cos t_\af)$, then an easy computation shows that the derivative of $a^{-p/2}I_\gamma|_{(b,x)=(1,0)}$ with respect to $a$ is positive at $\af$. This means that $\af$ is the {\it unique} value of $a$ that makes $I_\gamma|_{(b,x)=(1,0)}$ equals zero. Since the integral at (10) is increasing with $x$, we have just proved the following:
\\

{\it For any $a\in(0,\af)$, there exists a unique $x=x_a$ such that $I_\gamma|_{(b,x)=(1,x_a)}=0$. If $a\in(\af,1)$, then $I_\gamma|_{b=1}$ is always positive. Moreover, $\exists\Lim{a\to\af}{x_a=0}$.} 
\\

Let us now analyse the integral at (9). For $x=0$, it diverges to $+\infty$ when $a$ approaches 1. Take a compact $\K\subset\C\setminus\{0,a^{\pm 1}\}$ such that $S^1\subset\K$. One easily sees that our data $(1/g,dh)$ converge uniformly in $\K$ to the Weierstrass pair $(\bbg,\bbg\eta)$ from [MR,pp.452-3], for the following choice of parameters defined there: $\bba=1/a$, ${\boldsymbol A}=1$ and ${\boldsymbol B}=i\sqrt{\bba}$. Therefore, $I_\gamma|_{(b,x)=(a,0)}$ coincides with $\m\int_{\gamma_2}{{\bbf}_2}$, where $\gamma_2$ is described in [MR,pp.455]. There one proves that $\int_{\gamma_2}{{\bbf}_2}\ne 0$ for any $\bba>1$. Consequently, $I_\gamma|_{(b,x)=(a,0)}>0$ for all $a\in(0,1)$. Since $I_\gamma|_{b=a}$ is increasing with $x$, then $I_\gamma|_{b=a}>0$ on the whole square $(0,1)^2\ni(a,x)$.

We recall that $Re(g(t))$ is increasing and decreasing with $x$ and $b$, respectively. Hence, there is a function $b(a,x)$, defined in the region $\cR:=\{(a,x)\in(0,1)^2:x\le x_a\}$, such that $I_\gamma|_{b=b(a,x)}=0$ and non-zero elsewhere. Moreover, $b(a,x)$ can be continuously extended to $\dd\cR$ and $\Lim{(a,x)\to(0,0)}{b(a,x)=0}$. Henceforth in this section, the parameter $b$ will always represent {\it this} function.  

One easily sees that $\phi_2$ is purely imaginary for $-1<z<-x$ and $b<z<1$. From Figure 5(b) we get $Re\int_\delta\phi_2+Re\int_\gamma\phi_2=Re\int_\sg\phi_2$. From the above paragraph, the second integral is zero. Therefore, it remains to prove that either $Re\int_\delta\phi_2$ or $Re\int_\sg\phi_2$ vanishes for a subtable choice of $(a,x)\in\cR$. In order to accomplish this task, we shall make use of the following result:
\\
\\
{\bf Lemma 5.1.}\it \ The above defined $\af$ is bigger than $1/2$.
\\
\\
Proof.\rm \ Let us take $x=1$ at (10) and study the imaginary part of the function $z^3(1-az)/(z-a)$, for $z=e^{it}$, $0\le t\le\pi$. A simple reckoning shows that
\BE
   Im\biggl\{\frac{z^3(1/a-z)}{z/a-1}\biggl\}=
   \frac{\sin(2t)/a^2-2\sin(3t)/a+\sin(4t)}{(\cos t/a-1)^2+\sin^2t/a}.
\EE

If $a<1/2$, the derivative of (12) at either $t=0$ or $t=\pi$ is positive. Although it vanishes at both extremes for $a=1/2$, one still concludes that $Im\{\cdot\}$ is increasing there. For $a=1/2$, one rewrites the numerator of (12) as $4\sin t(1-\cos t)(2\cos t-\cos 2t)$. Since $\sin t(1-\cos t)$ never vanishes in $(0,\pi)$, in this interval there is a single zero at $t=t_0:=\arccos((1-\sqrt{3})/2)$. The real part of $z^3(1-az)/(z-a)$ has the same sign of $\cos(2t)/a^2-2\cos(3t)/a+\cos(4t)$, which at $t_0$ worths $-2^{\frac{3}{2}}\cdot 3^{\frac{1}{4}}\cdot(2\sqrt{3}-3)^\m-4\sqrt{3}$. This means that the argument of $z^3(1-az)/(z-a)$ varies from $0$ to $2\pi$ without taking negative values. Therefore, the integral at (10) is negative at $a=1/2$.\hfill q.e.d.  
\\

Now we parametrise the curve $\delta$ as $z(t)=z\circ\delta(t)=t,a<t<b$. If $\Delta:=\rho\circ\delta$ then $Z\circ\Delta(t)=z\circ\delta(t)$. From Figure 6 we have 
\BE
   g(\delta(t))=i|g(t)|=it^\frac{1}{4}\biggl(\frac{1-at}{t-a}\biggl)^\frac{1}{4}
   \biggl(\frac{b-t}{1-bt}\biggl)^\m\biggl(\frac{t+x}{xt+1}\biggl)^\m
\EE
and
\BE 
   dh(t)=|dh(t)|=\frac{dt/t}{\sqrt{a+1/a-t-1/t}}.
\EE

Therefore, $\phi_2\circ\delta(t)=(|g|^{-1}-|g|)|dh|$. On the points $(a,x)=(a,x_a)$ we have $b\equiv 1$. Under this condition and from (13), $\phi_2$ will be negative providing 
\BE
   a(t^4-1)+(x^2+2ax-1)(t^3-t)<0.
\EE

Since $t^2-1$ is always negative in $(a,1)$, then (15) is equivalent to
\BE
   t+1/t>(1-2ax-x^2)/a.
\EE

A sufficient condition for (16) to hold is that $2a>1-2ax-x^2$. Due to Lemma 5.1, it follows that $Re\int_\delta\phi_2$ is negative for $a$ close to $\af$. Now split $\sg$ into two stretches, the first one parametrised as $z(t)=z\circ\sg(t)=-t$, $0<t<x$, and the second $z(t)=t$, $0<t<a$. If $\Sg:=\rho\circ\sg$ then $Z\circ\Sg(t)=z\circ\sg(t)$. For the first stretch, from Figure 6 we have 
\BE
   g(\sg(t))=i|g(t)|=it^\frac{1}{4}\biggl(\frac{1+at}{a+t}\biggl)^\frac{1}{4}
   \biggl(\frac{b+t}{1+bt}\biggl)^\m\biggl(\frac{x-t}{1-xt}\biggl)^\m
\EE
and
\BE 
   dh(t)=|dh(t)|=\frac{dt/t}{\sqrt{t+1/t+a+1/a}}.
\EE
For the second stretch,
\BE
   g(\sg(t))=e^{\frac{i\pi}{4}}|g(t)|=
   e^{\frac{i\pi}{4}}t^\frac{1}{4}\biggl(\frac{1-at}{a-t}\biggl)^\frac{1}{4}
   \biggl(\frac{b-t}{1-bt}\biggl)^\m\biggl(\frac{x+t}{1+xt}\biggl)^\m
\EE
and
\BE 
   dh(t)=i|dh(t)|=\frac{dt/t}{\sqrt{t+1/t-a-1/a}}.
\EE

Thus $Re\int_\sg\phi_2=J_1-J_2$, where
\BE
   J_1:=\int_0^x\biggl(\frac{1}{|g|}-|g|\biggl)|dh|\eh\eh{\rm and}\eh\eh
   J_2:=\frac{\sqrt{2}}{2}\int_0^a\biggl(\frac{1}{|g|}+|g|\biggl)|dh|.
\EE

The change $t=au$ shows that $\Lim{a\to 0}{a^{-1/2}}J_2$ exists and is finite. Regarding $J_1$, from (17) we shall have $1/|g|>|g|$ providing $(b+t)(x-t)<(1+bt)(1-xt)$. This last inequality is equivalent to $t^2+2(b-x)t/(1-bx)+1>0$, which holds indeed, since $b-x>bx-1$. Now, an easy computation shows that $\Lim{a\to 0}{a^{-1/2}}J_1=+\infty$. We recall that $Re\int_\sg\phi_2=Re\int_\delta\phi_2$, and the latter is negative on $(a,x_a)$, $a$ close to $\af$. These facts imply that there is a curve $\cC\subset$ graph($b$) such that both $Re\int_\delta\phi_2$ and $Re\int_\gamma\phi_2$ vanish simultaneously for every choice of $(a,b,x)\in\cC$.
\\ 
\\
{\bf 7. Refinements} 
\\ 

In this section we study the curve $\cC$ with more details. First of all, let us prove 
\\
\\
{\bf Lemma 6.1.}\it \ There exists $\Lim{a\to 0}{x_a=1}$.
\\
\\
Proof.\rm \ From (10), an easy computation shows that
\BE
   \lim_{a\to 0}a^{-1/2}I_\gamma|_{b=1}=-\int_{0}^{\pi}
   Im\biggl\{\sqrt{\frac{z+x}{xz+1}}\biggl\}dt.
\EE

The integral at (22) is negative for any $x\in(0,1)$, but converges to zero when $x$ approaches 1. Since $\cR$ is exactly the region where $I_\gamma|_{b=1}$ is non-positive, the same holds for this integral re-scaled by $a^{-1/2}$. Suppose there were a positive $\eps$ admitting a sequence $a_n\to 0$ with $x(a_n)<1-\eps$ for all indexes $n$. In this case, the continuity of $a^{-1/2}I_\gamma|_{b=1}$, together with the fact that it is increasing with $x$, should give a non-negative limit in (22) for $x=1-\eps$. This would be a contradiction. Therefore, it exists $\Lim{a\to 0}{x_a=1}$.\hfill q.e.d.
\\

In the reminder of this section, we prove that the curve $\cC$ does not touch graph$(x_a)$. Hence, it connects the point $(0,1,1)$ with some point of graph$(b)$ over $(0,\af)\times\{0\}\ni(a,x)$. This will give a continuous one-parameter family of minimal surfaces with special limit-members. We shall describe them in Section 8. 

From (1), if we settle $b=1$, this gives another family of compact Riemann surfaces $R$ with algebraic equation
\BE    
   g^4=z\biggl(\frac{1-az}{z-a}\biggl)\biggl(\frac{z+x}{xz+1}\biggl)^2.  
\EE 
Of course, (23) {\it cannot} be viewed as a limit of (1) for $b\to 1$. The algebraic equations describe abstract surfaces, not even contained in a metric space. Our only resource is the study of period integrals, of which some limits can converge to integrals on another compact surface.   

The surfaces in (23) are endowed with the following involution: $(z,g)\to(z,ig)$. Since $i^4=1$, there are exactly four points of branch order 3, namely $(0,0)$, $(1/a,0)$, $(a,\infty)$ and $(\infty,\infty)$. Moreover, it remains only four other branch points, $(-x,\pm 0)$ and $(-1/x,\pm\infty)$, these of order 1. Here the $\pm$ signs indicate different {\it germs} of functions. The Riemann-Hurwitz formula gives
\[
   \frac{4\cdot 3+4\cdot 1}{2}-4+1=5.
\]

From now on, our analysis will be strongly based in [V1]. There one proves that the algebraic equations
\BE
   \biggl(G+\frac{1}{G}\biggl)^2=\frac{4\zeta(\zeta-y)^2(\zeta-1/\ld)}
   {(\zeta^2-1)(\zeta-\ka)(\zeta-1/\ka)}
\EE
and
\BE
   \biggl(G-\frac{1}{G}\biggl)^2=\frac{4(1-y\zeta)^2(1-\zeta/\ld)}
   {(\zeta^2-1)(\zeta-\ka)(\zeta-1/\ka)}
\EE 
are equivalent if and only if $\ld(\ka+1/\ka)=1+(2\ld-y)y$, with $2\ld-1<y<\ld<\ka<1$ and positive $\ld$. Moreover, the Riemann surfaces $M$ defined by (24-5) have genus 5. Notice that $M$ is endowed with the involution $\imath$ given by $(\zeta,G)\to(1/\zeta,iG)$. 

From [V1,pp351-3] one has that $\zeta$ is the pull-back under $\imath^2$ of an elliptic function $\cZ$ defined on a rectangular torus $\cT$. The parameter $\ld$ can freely vary in $(0,1)$, describing all rectangular tori. From Section 3 and [V1,pp352], one sees that the choice $a=\ld$ makes $T=\cT$ and $Z$ a ``shift'' of $\cZ$. By defining $\Ld:=\ld+1/\ld$, the following relation holds:
\BE
   \biggl(\frac{Z+1}{Z-1}\biggl)^2=\frac{\cZ+1/\cZ-\Ld}{2-\Ld}.
\EE 
If we choose $\cZ=\ka$, a unique $Z(\ka)\in(-1,0)$ will be determined by (26). So we take $a=\ld$ and $x=-Z(\ka)$ in (23). 

Let $\bbz$ be the pull-back of $Z$ under $\imath^2$. Therefore, $\bbz((1,0))=0$, $\bbz((1,\infty))=\infty$, $\bbz((-1,\infty))=a$ and $\bbz((-1,0))=1/a$, while $\bbz((\ka^{\pm 1},0))=-x$ and $\bbz((\ka^{\pm 1},\infty))=-1/x$. Let $\ell_j$ be a single small loop in $\C$ around $0$, $a$, $1/a$, $-x$ and $-1/x$, for $j=1,\dots,5$, respectively. We take lifts $\hat{\ell}_j$ of $\ell_j$ by $\bbz$ and notice that the end points of $\hat{\ell}_j$ differ by $\imath^{k_j}$, $0\le k_j\le 3$, $1\le j\le 5$.

Let $D$ be the open unitary complex disk at the origin. Since $deg(\bbz)=4$, there is a coordinate chart $w:D\to M$ with $w(0)=(1,0)$ such that $\bbz(w)=w^4$. By taking $\ell_1$ small enough to be in $\bbz(w(D))$, we conclude that $k_1=1$. The same reasoning will give $k_2=-k_3=-1$. If we had taken $w(0)=(\ka,0)$, then $\bbz(w)=c_1+w^2$ and so $k_4=2$. By the same reasoning $k_5=-2$. Let us define $\A:=\C\setminus\{0,a^{\pm 1},-x^{\pm 1}\}$.

The numbers $k_j$ naturally determine a homomorphism $H:\pi_1(\A)\to\Z_4\oplus\Z_2$, of which the kernel is $\bbz_*(\pi_1(M\setminus\{(\pm 1,0^{\pm 1}),(\ka^{\pm 1},0^{\pm 1})\}))\subset\pi_1(\A)$. By going back to (23), one sees that the projection map $z:R\to\hat{\C}$, namely $(z,g)\to z$, is such that $z_*(\pi_1(R\setminus g^{-1}(\{0,\infty\})))$ also represents the kernel of $H$. From [M,p159], there is a fibre-preserving biholomorphism $\bt:M\to R$ such that $\bbz=z\circ\bt$. As a matter of fact, that reference treats unbranched coverings, but the conclusion still applies to our case.     

From [V1] and the above paragraph, one sees that $G^4$ has the same divisor as $\bbz(1-a\bbz)/(\bbz-a)\cdot[(\bbz+x)/(x\bbz+1)]^2$. By composing $\bbz$ with the involution $(\zeta,G)\to(\bar{\zeta},1/\bar{G})$ we get $\bbz\to 1/\bar{\bbz}$. Therefore, $G$ is unitary where $\bbz$ is. Now, by composing $\bbz$ with the involution $(\zeta,G)\to(1/\bar{\zeta},i\bar{G})$ we get $\bbz\to\bar{\bbz}$. This means that $|\zeta|=1$ implies $\bbz^{\pm 1}\in(0,a)$. Hence
\[    
   G^4=\bbz\biggl(\frac{1-a\bbz}{\bbz-a}\biggl)\biggl(\frac{\bbz+x}{x\bbz+1}\biggl)^2,  
\] 
and so we can take $G=g\circ\bt$. In [V1] one defines $dH$ as the pull-back of the holomorphic differential form on $T$. As we have already mentioned, $\zeta$ is the pull-back of $\cZ$, which is a shift of $Z$. Hence, the pull-back of $\cZ'/\cZ$ gives a well-defined square-root of $\bbz+1/\bbz-a-1/a$ on $M$, and so $dH$ is proportional to $\bbz^{-1}d\bbz/\sqrt{\bbz+1/\bbz-a-1/a}$. But since $dH$ is purely imaginary for $|\bbz|=1$, the proportional constant must be $\pm i$. The sign just changes the minimal immersion to its antipodal, so we take
\[ 
   dH=\frac{id\bbz/\bbz}{\sqrt{\bbz+1/\bbz-a-1/a}}.
\] 

From Proposition 8.1 of [V1], or even better [V2], and the above discussion, one sees that each $a\in(0,1)$ admits a {\it unique} $\bbx$ for which $Re\int_{\bbsg}\bbf_2=0$. Now suppose that $\cC\cap$ graph$(x_a)\ne\emptyset$. In this case, there is $\ta\in(0,\af)$ such that $Re\int_\delta\phi_2=Re\int_\sg\phi_2=0$ for $(a,b,x)=(\ta,1,x_\ta)$. Therefore, it exists $\eps>0$ such that $Re\int_\delta\phi_2=0$ for all $(a,b,x)\in B_\eps(\ta,1,x_\ta)\cap\cC$.
 
Back to (13) and (14), the change $t\mapsto b-t^2$ shows that    
\[
   \exists\eh\lim_{b\to 1}Re\int_\delta\phi_2=Re\int_a^1\bbf_2.
\]
From the uniqueness, $\bbx=x_\ta$ because $Re\int_{\bbsg}\bbf_2=Re\int_\sg\phi_2|_{b=1}$ and the latter is zero at $(\ta,1,x_\ta)\in\cC$. But in [V1] one proves that such a choice gives an embedded surface, and in particular $Re\int_a^1\bbf_2$ is {\it negative}. 

Because of that, if $\eps$ is close enough to zero, then $Re\int_\delta\phi_2$ must be negative in $B_\eps(\ta,1,x_\ta)$, a contradiction. We conclude that $\cC\cap$ graph$(x_a)=\emptyset$. Consequently, the curve $\cC$ connects $(a,b,x)=(0,1,1)$ with $(a,b,x)=(a^*,b^*,0)$, for a certain $a^*\in(0,\af)$ and $b^*=b(a^*,0)$.       
\\ 
\\
{\bf 8. Limits and embeddedness} 
\\

At this point we have proved all but one item of Theorem 1.1. This last section is devoted to its accomplishment. By lopping off occasional loops of $\cC$, we can consider it as a simple curve. Let $s\mapsto(a(s),b(s),x(s))$ be a monotone parametrisation of $\cC$, assuming $(0,1,1)$ at $s=0$ and $(a^*,b^*,0)$ at $s=1$. For every $s\in(0,1)$, we have a well-defined minimal immersion $X_s:\bS\to\R^3/G$, determined by $(g,dh)$ at (1) and (4).

Now consider $u$ as a complex variable of $\hat{\C}$ and take $z=au/(u-1)$ in (1) and (4). If $K$ is a compact subset of $\hat{\C}\setminus\{1\}$, for $u\in K$ a simple computation gives
\BE
   \lim_{s\to 0}g^4=u\eh\eh{\rm and}\eh\eh
   \lim_{s\to 0}\frac{dh}{\sqrt{a}}=\frac{4dg/g}{g^2-1/g^2}.
\EE   
One readily recognises (27) as the Weierstrass data of Scherk's doubly periodic surface. Namely, the coordinates of the minimal immersion $X_s$ converge uniformly in $K$ to Scherk's coordinates. More precisely, suppose that $K$ is the 4$^{\rm th}$-power image of a compact $\sK\subset\hat{\C}\setminus\{\pm 1,\pm i\}$. In this set, $g$ is the standard complex coordinate, which together with $gdg/(g^4-1)$ gives the classical Scherk's doubly periodic surface. Figure 8 shows how the surface look like for $a$ close to zero.
\input epsf 
\begin{figure} [ht] 
\centerline{ 
\epsfxsize 11cm 
\epsfbox{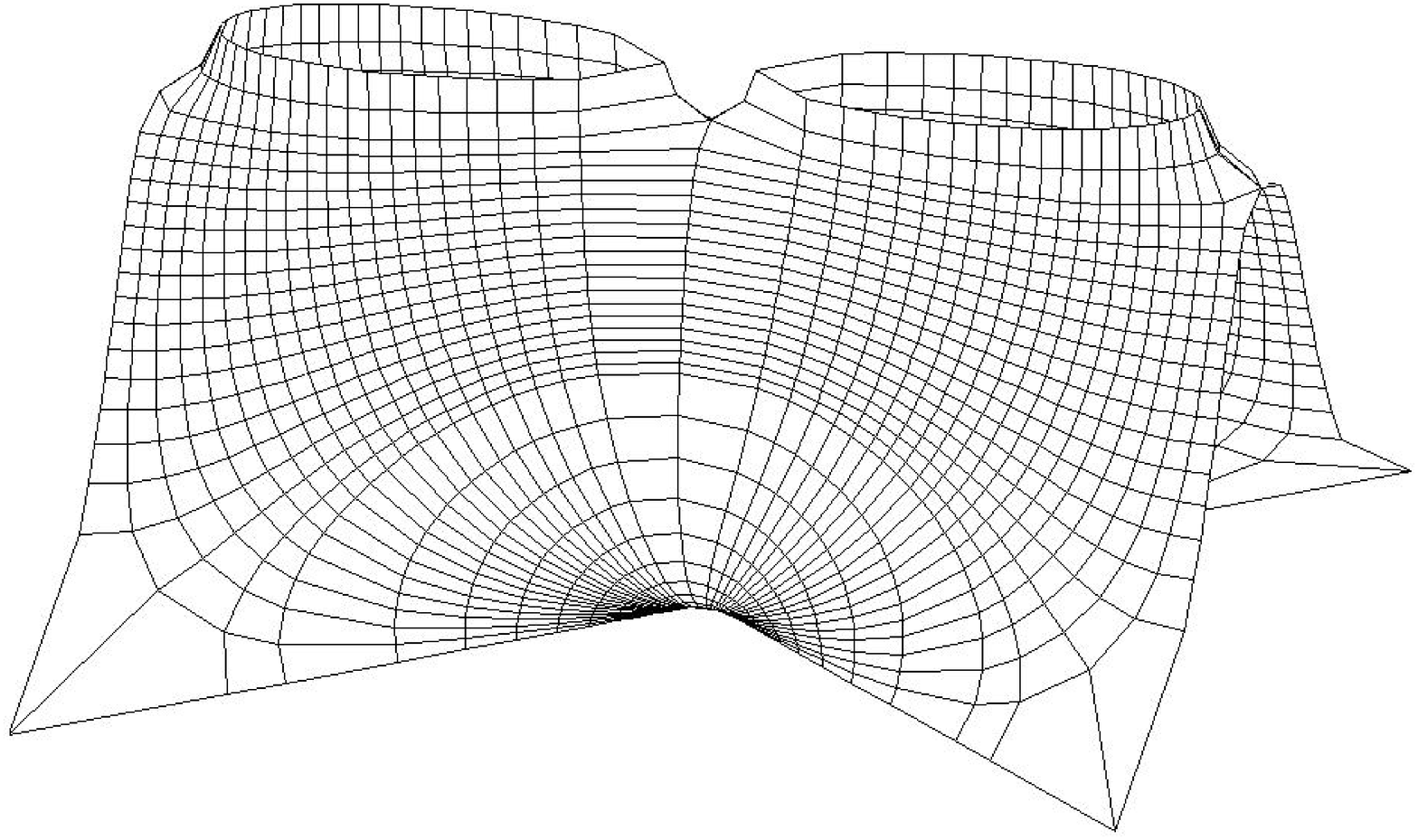}}
\caption{The case $(a,b,x)=(0.15;0.8;0.74)$.} 
\end{figure}

Consider now $z$ as complex variable of $\hat{\C}$ and define $\sD:=\{z\in\C:|z|<1<1+Im\{z\}\}$. For $z$ in a compact $\bbK\subset\sD\setminus\{0\}$, one immediately gets
\BE
   \lim_{s\to 1}g^4=z^3\biggl(\frac{1-a^*z}{z-a^*}\biggl)\biggl(\frac{b^*-z}{b^*z-1}\biggl)^2,
\EE 
while $\Lim{s\to 1}{dh}$ is given by (4) with $a=a^*$. From [SV,Sec.7] we recognise the Weierstrass data of a genus 5 example from Hoffman-Wohlgemuth. In fact, until the moment there is just numerical evidence that each genus 4$k$+1 gives a {\it unique} Hoffman-Wohlgemuth surface, $k\in\N^*$. However, in [SV] one gets all such surfaces from the {\it intermediate value theorem}. The choice $(a,b)=(a^*,b^*)$ is then included in [SV], since our surfaces are period free for all $s\in(0,1)$. 

Finally, the same arguments from [V1,p.360-2] imply that $X_s$ is in fact an embedding, for any $s\in(0,1)$. We conclude this last section with Figure 9, which illustrates the above convergence. Figure 1 shows the fundamental piece for $(a,b,x)=(0.47;0.85;0.68)$.
\input epsf 
\begin{figure} [ht] 
\centerline{ 
\epsfxsize 14cm 
\epsfbox{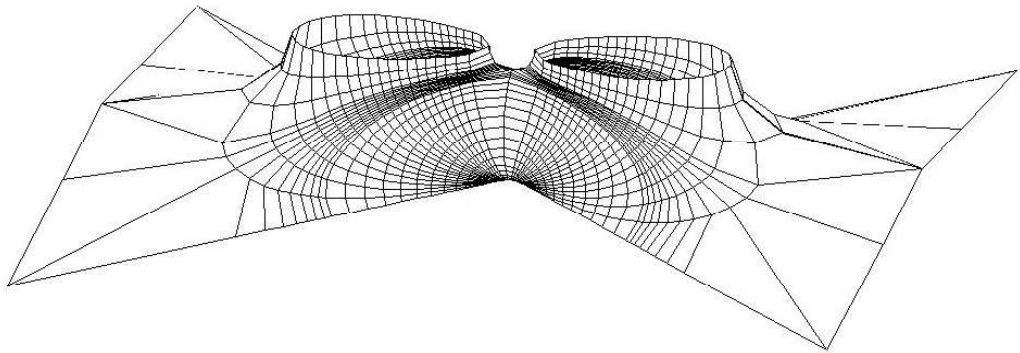}} 
\caption{The case $(a,b,x)=(0.65;0.89;0.69)$.} 
\end{figure} 
\ \\
{\bf References}
\ \\
\ \\
$[$C$]$- M. Callahan, D. Hoffman and W. H. Meeks. Embedded minimal surfaces with an infinite number of ends. {\it Inventiones Math.}, Vol.96, 1989, 459-505.
\ \\
$[$CK$]$- M. Callahan, D. Hoffman and H. Karcher. A family of singly periodic minimal surfaces invariant under a screw motion. {\it Experiment. Math.}, Vol.2, 1993, 157-182.
\ \\
$[$HK$]$- D. Hoffman, H. Karcher. {\it Complete embedded minimal surfaces of finite total curvature}, Encyclopedia of Math. Sci., Springer Verlag {\bf 90} (1997) 5--93.      
\ \\
$[$HKW$]$- Hoffman, David; Karcher, Hermann; Wei, Fu Sheng. The genus one helicoid and the minimal surfaces that led to its discovery. Global analysis in modern mathematics (1992), 119--170, Publish or Perish, Houston, TX, 1993. 
\ \\
$[$HPR$]$- L. Hauswirth, J. Perez and P. Romon. Embedded minimal ends of finite type. Trans. Amer. Math. Soc. 353 (2001), no. 4, 1335--1370.
\ \\
$[$HW$]$- D. Hoffman and H. Wohlgemuth. New embedded periodic minimal surfaces of Riemann-type. In manuscript, 1993.
\ \\
$[$K1$]$- H. Karcher. Embedded minimal surfaces derived from Scherk's examples. {\it Manuscripta Math.} {\bf 62} (1988), 83--114.
\ \\
$[$K2$]$- H. Karcher. The triply periodic minimal surfaces of Alan Schoen and their constant mean curvature companions. {\it Manuscripta Math.} {\bf 64} (1989), 291--357.
\ \\
$[$K3$]$- H. Karcher. Construction of minimal surfaces. {\it Surveys in Geometry}, University of Tokyo, 1989, 1-96 and {\it Lecture Notes}, Vol. 12, 1989, SFB256, Bonn.
\ \\
$[$LM$]$- F.J. L\'opez \& F. Mart\ih n, {\it Complete minimal surfaces in $\R^3$}, Publ. Mat. {\bf 43} (1999) 341--449.  
\ \\
$[$MR$]$- F. Mart\ih n \& D. Rodr\ih guez. A characterization of the periodic Callahan-Hoffman-Meeks surfaces in terms of their symmetries. {\it Duke Math. J.}, Vol. 89, 1997, 445-463.
\ \\
$[$N$]$- J.C.C. Nitsche, {\it Lectures on minimal surfaces}, Cambridge University Press, Cambridge (1989).              
\ \\
$[$O$]$- R. Osserman, {\it A survey of minimal surfaces}, Dover, New York, 2nd ed (1986).
\ \\
$[$PRT$]$- J. P\'erez, M. Rodr\ih guez \& M. Traizet, {\it The classification of doubly periodic minimal tori with parallel ends}, J. Differential Geom. {\bf 69} (2005) 523--577.
\ \\
$[$SV$]$- P.A.Q.Sim\~oes and R.B. Val\'erio. A characterisation of the Hoffman-Wohlgemuth surfaces in terms of their symmetries. J. Differential Geom., to appear.
\ \\
$[$V1$]$- R.B. Val\'erio. A family of triply periodic Costa surfaces. {\it Pacific J. Math.}, Vol.212, 2003, 347-370.
\ \\
$[$V2$]$- R.B. Val\'erio. Theoretical evaluation of elliptic integrals based on computer graphics. UNICAMP Technical Report 71/02, Campinas, SP 2002; home page http://www.ime.unicamp.br/rel\_pesq/2002/rp71-02.html
\ \\
$[$HWW$]$- D. Hoffman, M. Wolf and M. Weber. The genus-one helicoid as a limit of screw-motion invariant helicoids with handles. Clay Math. Proc., 2, Global theory of minimal surfaces,  243--258, Amer. Math. Soc., Providence, RI, 2005. 
\ \\
$[$W$]$- F. Wei. Some existence and uniqueness theorems for doubly periodic minimal surfaces. {\it Inventiones Math.}, Vol.109, 1992, 113-136.
\ \\
$[$M$]$- W.S. Massey. Algebraic topology: an introduction. {\it Graduate Texts in Mathematics}, Springer, New York (1967).
\end{document}